\documentclass[twoside,a4paper]{article}

\usepackage{amsmath,amsthm,amssymb,latexsym}
\usepackage[v2,tips]{xy}

\newcommand{\wap}{\operatorname{WAP}}
\newcommand{\ip}[2]{{\langle {#1} , {#2} \rangle}}

\newcommand{\mc}[1]{\mathcal{#1}}

\newcommand{\aone}{\Box}
\newcommand{\atwo}{\Diamond}

\theoremstyle{plain}%
\newtheorem{proposition}{Proposition}[section]%
\newtheorem{theorem}[proposition]{Theorem}%
\newtheorem{corollary}[proposition]{Corollary}%
\newtheorem{lemma}[proposition]{Lemma}%
\theoremstyle{definition}%
\newtheorem{definition}[proposition]{Definition}%
\theoremstyle{remark}%
\newtheorem{remark}[proposition]{Remark}%

\begin{document}

\large
\title{\textsc{Preduals of semigroup algebras}}
\author{Matthew Daws\\\normalsize{\texttt{matt.daws@cantab.net}}\and Hung Le Pham\\\normalsize{\texttt{hlpham@math.ualberta.ca}}\and Stuart White\\\normalsize{\texttt{s.white@maths.gla.ac.uk}}}
\maketitle

\begin{abstract}
For a locally compact group $G$, the measure convolution algebra $M(G)$ carries a natural coproduct.  In previous work, we showed that the canonical predual $C_0(G)$ of $M(G)$ is the unique predual which makes both the product and the coproduct on $M(G)$ weak$^*$-continuous.  Given a discrete semigroup $S$, the convolution algebra $\ell^1(S)$ also carries a coproduct.  In this paper we examine preduals for $\ell^1(S)$ making both the product and the coproduct weak$^*$-continuous. Under certain conditions on $S$, we show that $\ell^1(S)$ has a unique such predual. Such $S$ include the free semigroup on finitely many generators.  In general, however, this need not be the case even for quite simple semigroups and we construct uncountably many such preduals on $\ell^1(S)$ when $S$ is either $\mathbb Z_+\times\mathbb Z$ or $(\mathbb N,\cdot)$.

\end{abstract}

\section{Introduction}

A \emph{dual Banach algebra} is a Banach algebra $\mc A$ which is the dual of a
Banach space $\mc A_*$, such that the product on $\mc A$ is separately
weak$^*$-continuous.  The motivating example is a von Neumann algebra, where the
predual is isometrically unique.  This need not be true for Banach algebras:
consider $\ell^1$ with the zero product.  In \cite{our}, we considered the measure
algebra of a locally compact group, $M(G)$, which is a dual Banach algebra with
respect to the predual $C_0(G)$.  We can define a natural coproduct on $M(G)$,
dualising the multiplication on $C_0(G)$.  In \cite[Theorem~3.6]{our}, we showed
that $C_0(G)$ is the unique predual making both the product and the coproduct on
$M(G)$ weak$^*$-continuous.  The proof makes use of results which are only true
for group algebras.  In this paper we consider preduals of $\ell^1$-semigroup algebras
which make both the product and coproduct weak$^*$-continuous. Such preduals are termed
\emph{Hopf algebra preduals}.  Perhaps surprisingly, we show that quite simple
semigroups, such as $\mathbb Z_+\times\mathbb Z$, give rise to algebras with
uncountably many Hopf algebra preduals.

Given a predual $\mc A_*$ of $\mc A$, we can naturally identify $\mc A_*$ with a
closed subspace of $\mc A^*$.  Then the product on $\mc A$ becomes separately
weak$^*$-continuous if and only if $\mc A_*$ is a submodule of $\mc A^*$, for the
usual action of $\mc A$ on its dual.  It is easy to see that isomorphic preduals
will induce the same subspace of $\mc A^*$.  Henceforth, by a \emph{predual}
$\mc A_*$ of a Banach algebra $\mc A$, we shall mean a closed submodule of
$\mc A^*$ which is a Banach space predual.  

We shall consider convolution algebras $\ell^1(S)$ for a (countable) semigroup $S$.
We write $(\delta_s)_{s\in S}$ for the standard unit vector basis of $\ell^1(S)$,
so each $a\in\ell^1(S)$ can be uniquely expressed as a norm-convergent sum
\[ a = \sum_{s\in S} a_s \delta_s \quad\text{where}\quad
\|a\| = \sum_{s\in S} |a_s|. \]
The coproduct on $\ell^1(S)$ is the map $\Gamma:\ell^1(S)\rightarrow\ell^1(S\times S)$
defined by
\[ \Gamma(\delta_s) = \delta_{(s,s)} \qquad (s\in G). \]
Let $E\subseteq\ell^\infty(S) = \ell^1(S)^*$ be a predual for $\ell^1(S)$.  Then
it is easy to see (compare with \cite[Lemma~3.3]{our}) that $\Gamma$ is weak$^*$-continuous
if and only if $E$ is a subalgebra of $\ell^\infty(S)$.  As in \cite{our}, we shall term
such a predual a \emph{Hopf algebra predual} of $\ell^1(S)$.

In the next section, we show that the study of Hopf algebra preduals of $\ell^1(S)$ is equivalent to the study of certain semigroup topologies on $S$. For certain cancellative semigroups $S$, including a finite direct sum of copies of $\mathbb Z_+$ and the free semigroup on finitely many generators, we show that $\ell^1(S)$ has a unique Hopf algebra predual.  In Section~\ref{count}, we exhibit semigroups admitting uncountably many Hopf algebra preduals.  These semigroups are quite simple in nature, including $\mathbb Z_+\times\mathbb Z$ and $(\mathbb N,\cdot)$.  Since this last semigroup is isomorphic to the direct sum of countably many copies of $\mathbb Z_+$, it is not possible to extend the uniqueness result for Hopf algebra preduals of finite direct sums of $\mathbb Z_+$ to infinite direct sums.  In Section~4, we exhibit a semigroup $S$ for which $\ell^1(S)$ admits no Hopf algebra preduals and we end in Section~5 by showing that $\ell^1(\mathbb N,\max)$ has a unique predual in full generality. 

Finally, some words about notation.  For a Banach space $E$, we write $E^*$ for its
dual, and use the dual-pairing notation $\ip{\mu}{x} = \mu(x)$ for $\mu\in E^*$ and
$x\in E$.

\vspace{0.2cm}
\noindent\textbf{Acknowledgments.}
The second named author is supported by a Killam Postdoctoral Fellowship and a
Honorary PIMS PDF.  He also wishes to thank Professor Anthony To-Ming Lau for his
kind support and encouragement during this research.

\section{General theory}

Let us recall the following, which is \cite[Lemma~3.4]{our}.

\begin{lemma}\label{cstar_bidual_is_cstar}
Let $K$ be a locally compact Hausdorff space, let $L$ be a compact Hausdorff
space, and let $\mc A$ be a closed subalgebra of $C_0(K)$ such that $\mc A^{**}$
is Banach algebra isomorphic to $C(L)$.  Then $\mc A$ is a C$^*$-subalgebra of $C_0(K)$.
\end{lemma}

A semigroup $K$ which carries a topology is said to be \emph{semitopological}
when the product is separately continuous.

\begin{proposition}\label{reduce_to_topology}
Let $S$ be a discrete semigroup, and let $E\subseteq\ell^\infty(S)$ be a Hopf
algebra predual for $\ell^1(S)$.
Then $E$ is a C$^*$-subalgebra of $\ell^\infty(S)$.  Let $K$ be the character
space of $E$.  Then $K$ is canonically bijective to $S$, and we use this bijection
to turn $K$ into a semigroup.  Then:
\begin{enumerate}
\item $K$ is a semitopological semigroup;
\item\label{rtt:sec} for each $s\in K$ and $L\subseteq K$ compact, the sets
   $\{ t\in K : st\in L\}$ and $\{ t\in K : ts\in L\}$ are compact.
\end{enumerate}
This structure on $K$ completely determines $E$.
\end{proposition}
\begin{proof}
Arguing as in the proof of \cite[Corollary~3.5]{our}, $E^{**}$ is Banach algebra
isomorphic to $\ell^\infty(S)$.  It then follows immediately from
Lemma~\ref{cstar_bidual_is_cstar} that $E$ is a C$^*$-subalgebra of $\ell^\infty(S)$.
Let $K$ be the character space, and let $\mc G:E\rightarrow C_0(K)$ be the Gelfand
transform.  Then following the proof of \cite[Theorem~3.6]{our},
we get a natural bijection $\theta:S\rightarrow K$ such that
\[ f(k) = \ip{\mc G^{-1}(f)}{\delta_{\theta^{-1}(k)}}
\qquad (f\in C_0(K), k\in K). \]
It also follows that $K$ is semitopological.

The second condition follows as the left and right translations by elements of
$K$ will take $C_0(K)$ to $C_0(K)$.  Let $L\subseteq K$ be compact, and
let $s\in K$.  For each $r\not\in L$, we can find $f_r\in C_0(K)$ with
$f_r(t)=1$ for $t\in L$, and $f_r(r)=0$.  Then let $g = f_r\cdot\delta_s$, so that
\[ g(t) = \ip{\delta_t}{g} = \ip{\delta_t}{f_r\cdot\delta_s}
= \ip{\delta_{st}}{f_r} = f_r(st) \qquad (t\in K). \]
Then as $E = C_0(K)$ is an $\ell^1(S)$-submodule, we have that $g\in C_0(K)$,
so the set $L_r = \{ t\in K : |f_r(st)|\geq1 \}$ is compact.  As $f_r(r)=0$, we
see that $t\in L_r$ implies that $st\not=r$.  Similarly, if $st\in L$, then
certainly $t\in L_r$.  So we conclude that
\[ \bigcap_{r\not\in L} L_r = \{t\in K : st\in L\} \]
is a compact set, as required.  Similarly, $\{t\in K : ts\in L\}$ is compact.

Finally, suppose that $F\subseteq\ell^\infty(S)$ is another predual, inducing a
semitopological semigroup $L$ by a bijection $\phi:S \rightarrow L$.  Suppose that $K=L$
in the sense that $\phi \theta^{-1}:K\rightarrow L$ is a bi-continuous homomorphism.
Let $x\in E$, so that
\[ \ip{x}{\delta_s} = \mc G(x)(\theta(s)) \qquad (s\in S). \]
Let $f = \mc G(x) \in C_0(K)$, and let $g \in C_0(L)$ be defined by
\[ g(l) = f\big( \theta\phi^{-1}(l) \big) \qquad (l\in L). \]
Let $\mc H:F\rightarrow C_0(L)$ be the Gelfand transform, and let $y=\mc H^{-1}(g)
\in F$.  Then
\[ \ip{y}{\delta_s} = \mc H(y)(\phi(s)) = g(\phi(s)) = f(\theta(s)) = \ip{x}{\delta_s}, \]
so that $x = y$.  As $x\in E$ was arbitrary, we conclude that $E\subseteq F$, and
by symmetry, that actually $E=F$, as required.
\end{proof}

\begin{lemma}
Let $S$ be a countable semigroup, and let $\sigma$ be a locally compact, semitopological
topology on $S$ satisfying condition (\ref{rtt:sec}) above.  Then $C_0(S,\sigma)$
is a predual for $\ell^1(S)$.
\end{lemma}
\begin{proof}
As $S$ is countable, it follows that $C_0(S,\sigma)^* = M(S,\sigma) = \ell^1(S)$.
Let $f\in C_0(S,\sigma)$, so that
\[ \ip{\delta_s\cdot f}{\delta_t} = f(ts) \qquad (s,t\in S), \]
so $\delta_s\cdot f\in C_0(S,\sigma)$ if and only if $f_s:t\mapsto f(ts)$ is
in $C_0(S,\sigma)$.  As $S$ is semitopological, $f_s$ is continuous.  For
$\epsilon>0$, we see that
\[ X = \{ s\in S : |f(s)|\geq\epsilon \} \]
is compact, so that
\[ \{ t\in S : |f_s(t)|\geq\epsilon \} = \{ t\in S : ts\in X \}, \]
is compact, by condition (\ref{rtt:sec}).  So $f_s\in C_0(S,\sigma)$.
Similarly, $f\cdot\delta_s\in C_0(S,\sigma)$ for $s\in S$, so $C_0(S,\sigma)$
is a submodule of $\ell^\infty(S)$.
\end{proof}

It hence follows that studying Hopf algebra preduals of $\ell^1(S)$ is completely
equivalent to studying locally compact semitopological topologies on $S$ which satisfy
condition (\ref{rtt:sec}), a task we shall concern ourselves with in
much of the rest of this paper.

\begin{corollary}\label{good_corr}
Let $S$ be a countable semigroup.  Then a locally compact topology $\sigma$ on $S$
makes $C_0(S,\sigma)$ into a predual for $\ell^1(S)$ if and only if:
\begin{enumerate}
\item $(S,\sigma)$ is a semitopological semigroup;
\item is $s\in S$ and $(t_n)\subseteq S$ is a sequence such that either $(t_ns)$ or
   $(st_n)$ has a convergent subsequence, then $(t_n)$ has a convergent subsequence.
\end{enumerate}
\end{corollary}
\begin{proof}
As $S$ is countable, $(S,\sigma)$ is metrisable.  Hence the second condition is
easily seen to be equivalent to the second condition in the result above.
\end{proof}

We say that a semigroup $S$ is \emph{weakly cancellative} when the left and right translation maps
are all finite-to-one maps; $S$ is \emph{cancellative} when these maps are all injective.
A simple calculation (see, for example, \cite[Theorem~4.6]{DL})
shows that $S$ is weakly cancellative if and only if $c_0(S)$ is a predual for $\ell^1(S)$.

However, even when $S$ is not weakly-cancellative, $\ell^1(S)$ might still have a predual: for
example, the one-point compactification of $\mathbb N$, denoted by $\mathbb N_\infty$, carries
an obvious semigroup structure for which $c$, the space of convergence sequences, forms a
predual for $\ell^1(\mathbb N_\infty)$.  This fits into our framework, as
$c = C_0(\mathbb N_\infty,\sigma)$ where $\sigma$ is the one-point compactification
topology on $\mathbb N_\infty$.

\begin{proposition}
Let $E\subseteq\ell^\infty(\mathbb N_\infty)$ be a Hopf algebra predual for
$\ell^1(\mathbb N_\infty)$.  Then $E = c$.
\end{proposition}
\begin{proof}
Let $(\mathbb N_\infty,\sigma)$ be the character space of $E$, so $\sigma$
satisfies the conditions of Corollary~\ref{good_corr}.  To show that $E=c$, we
need to show that $\sigma$ is the one-point compactification topology.
That is, we need to show that $\infty$ is the only limit point of $\sigma$.

Towards a contradiction, suppose that $(t_n)$ is a sequence in $\mathbb N$,
converging to $t<\infty$.  By moving to a subsequence, we may suppose that
$(t_n)$ is an increasing sequence, and that $t_n>t$ for all $n$.  Then
$(t_n) = ( t_n-t + t )$ converges, so by the corollary, $(t_n-t)$ have a
convergent subsequence, say $(s_n)$.  Let $s = \lim_n s_n$, so that as
$\sigma$ is semitopological,
\[ s+t = \lim_n s_n + t = \lim_n (t_n-t)+t = \lim_n t_n = t. \]
However, this contradicts $s\in\mathbb N = \{1,2,3,\cdots\}$.

\end{proof}

\subsection{Cancellation properties}

We shall now consider cancellative semigroups $S$, so that $c_0(S)$ is a
Hopf algebra predual for $\ell^1(S)$.

\begin{lemma}\label{best_semigroup}
Let $S$ be a countable cancellative semigroup.  Then a locally compact topology
$\sigma$ on $S$ makes $C_0(S,\sigma)$ into a predual for $\ell^1(S)$ if and only if:
\begin{enumerate}
\item $(S,\sigma)$ is a semitopological semigroup;
\item is $s\in K$ and $(t_n)\subseteq S$ is a sequence such that either $(t_ns)$ or
   $(st_n)$ converges, then $(t_n)$ converges.
\end{enumerate}
\end{lemma}
\begin{proof}
We need only check that this new second condition is implied by the second condition
in Corollary~\ref{good_corr}.  So let $s\in S$ and let $(t_n)\subseteq S$ be a sequence such
that $(t_ns)$ is convergent.  Hence $(t_n)$ has a convergent subsequence, converging
to $t_0\in S$ say.  As $(S,\sigma)$ is a semitopological semigroup, $(t_ns)$ converges
to $t_0s$.  Suppose that $(t_n)$ does not converge to $t_0$, so we can find a subsequence
which is always distant from $t_0$.  By applying the second condition from the
preceding lemma again, we find a different limit point, say $t_1$, for some further
subsequence.  Again, we hence have that $(t_ns)$ converges to $t_1s$.  So $t_0s=t_1s$,
and as $S$ is cancellative, $t_0=t_1$, a contradiction.  Hence $(t_n)$ converges,
as required.
\end{proof}

We make the following temporary definition.

\begin{definition}
Let $S$ be a semigroup.  We say that $S$ is \emph{finitely left divisible}
if for each $s\in S$, the set $\{ t\in S : \exists\,r\in S, tr=s \}$ is finite.
\end{definition}

\begin{theorem}
Let $S$ be a finitely left divisible, cancellative, finitely generated semigroup.
Let $E\subseteq\ell^\infty(S)$ be a Hopf algebra predual for $\ell^1(S)$.
Then $E = c_0(S)$.
\end{theorem}
\begin{proof}
By the previous results, we need to show that if $\sigma$ is a locally compact
topology on $S$ satisfying the conditions of Lemma~\ref{best_semigroup},
then $\sigma$ is the discrete topology.  Towards a contradiction, suppose that
$\sigma$ is not the discrete topology.  Hence we can find a non-isolated point
$t_0\in S$.  Let $(t_n)$ be some sequence in $S\setminus\{t_0\}$ which converges to $t_0$.

Let $C$ be a finite set which generates $S$.  For each $t\in S$, define the
\emph{length} of $t$ to be $l(t)$, the smallest $n$ such that $t=c_1\cdots c_n$
for $(c_i)\subseteq C$.  As $\{t_n\}$ is infinite, we see that $\{ l(t_n) \}$
is unbounded.  By moving to a subsequence, we may suppose that $l(t_1) < l(t_2)
< \cdots$.  Write $t_n = c_{n,1} c_{n,2} \cdots c_{n,l(t_n)}$ where $c_{n,i}\in
C$ for all $n$ and $i$.  By moving to a subsequence, we may suppose that
$c_{n,1} = c_1$ for all $n$.  By moving to a further subsequence, we may suppose
that $c_{n,2} = c_2$ for all $n$.  We can continue, and then by a diagonal argument,
we may suppose that $c_{n,k} = c_k$ for all $n$ and all $k\leq n$.

For each $m$, for $n\geq m$, we have that $t_n = c_1 c_2 \cdots c_n r_n$ for
some $r_n\in S$.  Let $u_{n,k} = c_k c_{k+1} \cdots c_n r_n$.
As $t_n \rightarrow t_0$, we see that $\lim_n c_1c_2\cdots c_{k-1} u_{n,k} = t_0$,
and so the corollary implies that $(u_{n,k})_{n=1}^\infty$ converges, to $v_k$ say.
Hence $t_0 = c_1c_2\cdots c_{k-1} v_k$, for all $k$.

As $S$ is finitely left divisible, we see that the set
$\{ c_1c_2\cdots c_k : k\geq 1 \}$
is a finite set in $S$.  In particular, we can find $j<k$ with
\[ c_1c_2\cdots c_j = c_1c_2\cdots c_k. \]
We have that $t_k = c_1c_2\cdots c_k c_{k,k+1}\cdots c_{k,l(t_k)}$,
where $l(t_k)$ is minimal.  However, we now also see that
$t_k = c_1c_2\cdots c_j c_{k,k+1}\cdots c_{k,l(t_k)}$, which is a shorter
expression, contradicting the minimality of $l(t_k)$.  This contradiction shows
that $t_0$ is indeed an isolated point, as required.
\end{proof}

\begin{corollary}\label{cor:one}
Let $k\geq 1$, and let $S=\mathbb Z_+^k$, or $S=\mathbb S_k$, the free semigroup
on $k$ generators.  Then $c_0(S)$ is the unique Hopf algebra predual of $\ell^1(S)$.
\end{corollary}

See Section~\ref{count} below for counter-examples showing that we cannot remove
the ``finitely left divisible'' or ``finitely generated'' conditions from the above
theorem.

\subsection{Rees semigroups of matrix type}

Rees matrix semigroups appear naturally in the study of when $\ell^1$-semigroup
algebras are amenable: see \cite[Chapter~10]{DL} for further details.  In general,
they are a way of generating semigroups from other (semi)groups.

Let $S$ be a semigroup with $0$, let $I$ and $J$ be index sets, let
$P$ be an $J\times I$ matrix of entries from $S$, called the \emph{sandwich
matrix}.  Then the \emph{Rees semigroup} $\mc M(S;I,J;P)$ is the collection
of all $I\times J$ matrices with entries from $S$, where exactly one entry is
non-zero.  Then we can define a product on $\mc M(S;I,J;P)$ by
\[ A\cdot B = APB \qquad (A,B\in\mc M(S;I,J;P)). \]
Of course, we cannot multiply matrices with entries in $S$, as we have
no concept of addition.  However, a moment's though reveals that in all
calculations, at most one entry will be non-zero, so there is no ambiguity
in what we mean by ``matrix multiplication'' in this setting.

An alternative description is the following.  Let
$(\epsilon_{i,j})_{i\in I, j\in J}$ be the matrix units, so each element of
$\mc M(S;I,J;P)$ can be written, formally, as $s \epsilon_{i,j}$ for some
$s\in S$.  Then
\[ s\epsilon_{i,j} \cdot t\epsilon_{k,l} = s P_{j,k} t \epsilon_{i,l}
\qquad (s,t\in S, i,k\in I, j,l\in J). \]


\begin{lemma}
For a group $G$, we have that $\mc M(G;I,J;P)$ is weakly-cancellative if and
only if $I$ and $J$ are finite.
\end{lemma}
\begin{proof}
If $J$ is infinite, then for $s\in G$, $i\in I$ and $j\in J$, fixed, notice that
\[ s \epsilon_{i,j} \cdot t \epsilon_{k,l}
= s P_{j,k} t \epsilon_{i,l}
= s P_{j,a} P_{j,a}^{-1} P_{j,k} t \epsilon_{i,l}
= s \epsilon_{i,j} \cdot P_{j,a}^{-1} P_{j,k} t \epsilon_{a,l}, \]
for any $t\in G$, $k,a\in I$ and $l\in J$.  So the map given by multiplication
on the left by $s\epsilon_{i,j}$ is an infinite-to-one map.  Hence $\mc M(G;I,J;P)$
weakly-cancellative implies that $J$ is finite; similarly $I$ must be finite.
\end{proof}

\begin{proposition}\label{Rees_prop}
Let $G$ be a countable group, let $I$ and $J$ be finite, and let
$S=\mc M(G;I,J;P)$.  Then $c_0(S)$ is the unique C$^*$-predual for $\ell^1(S)$.
\end{proposition}
\begin{proof}
Let $\sigma$ be a locally compact Hausdorff topology on $S$, making $(S,\sigma)$
a semitopological semigroup.  As $(S,\sigma)$ is a countable Baire space,
there exists $s_0 = g_0\epsilon_{i_0,j_0}\in S$ with $\{ s_0 \}$ open.
Then, for $g\in G$, $i\in I$ and $j\in J$, notice that
\begin{align*} \{ (h,k,l)\in G\times I\times J :
   h \epsilon_{k,l} \cdot g \epsilon_{i,j} = s_0 \}
&= \{ (h,i_0,l) : l\in J, h P_{l,i} g = g_0 \} \\
&= \{ (g_0 g^{-1} P_{l,i}^{-1}, i_0, l) : l\in J \} \end{align*}
if $j=j_0$, or the empty set otherwise.  Hence, for all $g\in G$ and $i\in I$,
the set
\[ U_{g,i} = \{ g P_{l,i}^{-1} \epsilon_{i_0, l} : l\in J \} \]
is open.  Similarly, we can check that
\[ V_{g,j} = \{ P_{j,k}^{-1} g \epsilon_{k,j_0} : k\in I \} \]
is open, for all $g\in G$ and $j\in J$.

Notice now that for $g,h\in G$, $I\in I$ and $j\in J$,
\[ U_{g,i} \cap V_{h,j} = \begin{cases} \big\{ g P_{j_0,i}^{-1}
  \epsilon_{i_0,j_0} \big\} &: g P_{j_0,i}^{-1} = P_{j,i_0}^{-1} h \\
  \emptyset &: \text{otherwise}. \end{cases} \]
So we have that $\{ g \epsilon_{i_0,j_0} \}$ is open, for any $g\in G$.

Let $X$ be the collection of open singletons in $(S,\sigma)$.  Then we have
just proved that if $g_0 \epsilon_{i_0,j_0}\in X$ for some $g_0$, then
$g \epsilon_{i_0,j_0}\in X$ for all $g\in G$.  Also, $X$ must be dense in
$(S,\sigma)$, for otherwise, the complement of the closure of $X$ would be a
non-empty open subset of a locally compact space, and hence locally compact
itself.  Repeating the Baire space argument would then yield an open singleton
not in $X$, a contradiction.

Towards a contradiction, suppose that $X$ is not all of $S$.  Then there
exists $s_0=g_0\epsilon_{i_0,j_0}\not\in X$ and a sequence
$(g_n\epsilon_{i_n,j_n})$ in $X$ converging to $s_0$.  By moving to a
subsequence, we may suppose that $i_n=i$ and $j_n=j$ for all $n$, so that
\[ g_n \epsilon_{i,j} \rightarrow s_0. \]
We then see that for $g,h\in G$, $k,p\in I$ and $l,q\in J$,
\begin{align*} g P_{l,i} g_n P_{j,p} h \epsilon_{k,q} 
&= g \epsilon_{k,l} \cdot g_n \epsilon_{i,j} \cdot h \epsilon_{p,q}
\rightarrow
g \epsilon_{k,l} \cdot g_0 \epsilon_{i_0,j_0} \cdot h \epsilon_{p,q}
= g P_{l,i_0} g_0 P_{j_0,p} h \epsilon_{k,q}. \end{align*}

Pick $(k,q)\in I\times J$ and $f\in G$ such that $f \epsilon_{k,q}
\in X$.  By the above, we see that
\[ g_0^{-1} P_{l,i_0}^{-1} P_{l,i} g_n P_{j,p} P_{j_0,p}^{-1} f \epsilon_{k,q}
\rightarrow
g_0^{-1} P_{l,i_0}^{-1} P_{l,i_0} g_0 P_{j_0,p} P_{j_0,p}^{-1} f \epsilon_{k,q}
= f \epsilon_{k,q}. \]
As $\{f \epsilon_{k,q}\}$ is open, we must have that
\[ P_{l,i} g_n P_{j,p} = P_{l,i_0} g_0 P_{j_0,p} \qquad (n\in\mathbb N). \]
However, this would imply that the
\[ P_{l,i}^{-1} P_{l,i_0} g_0 P_{j_0,p} P_{j,p} \epsilon_{i,j} \rightarrow
s_0 = g_0 \epsilon_{i_0, j_0}, \]
which is a contradiction, as $P_{l,i}^{-1} P_{l,i_0} g_0 P_{j_0,p} P_{j,p}
\epsilon_{i,j}\in X$, yet $s_0 \not\in X$.
\end{proof}

\section{Semigroup algebras with many preduals}\label{count}

In this section, we shall show that if $S=(\mathbb N,\cdot)$, the semigroup of natural
numbers with multiplication product, or if $S=(\mathbb Z_+ \times \mathbb Z,+)$, then
$\ell^1(S)$ admits a continuum of Hopf algebra preduals.  Notice that both semigroups are
cancellative, and the former is finitely (left) divisible but not finitely generated,
whereas the latter is finitely generated but not finitely (left) divisible. We begin with a elementary technical observation. 

\begin{lemma}\label{technical}
Let $a,\alpha\in\mathbb N$, and define
\[ X_{a,\alpha} = \{0\}\ \cup\ \Big\{\sum_{i=1}^k 2^{-m_i}:\ 1\le
k\le a,\ \alpha\le m_1<\ldots < m_k \Big\} \subseteq [0,1], \]
equipped with the subspace topology.  Let $\alpha\leq m_1 < m_2 < \cdots < m_k$ for some $1\leq k\leq a$, and
let $x_0 = \sum_{i=1}^k 2^{-m_i}$.  For $\beta> m_k$, let
\[ Y = \Big\{ \sum_{i=1}^{k+l} 2^{-m_i} : 0\leq l\leq a-k, \beta \leq m_{k+1} < \cdots < m_{k+l} \Big\}. \]
Then $Y$ is open in $X_{a,\alpha}$.
\end{lemma}
\begin{proof}
We claim that $(x_0-2^{-m_k-a},x_0]\cap X_{a,\alpha}=\{x_0\}$.  Given $x_1\in X_{a,\alpha}$ with $x_0-2^{-m_k-a}<x_1\leq x_0$, write 
$$
x_1=\sum_{i=1}^l 2^{-n_i}
$$
for some $1\leq l\leq a$ and $\alpha\leq n_1<\cdots<n_l$.  Certainly $n_1\geq m_1$, as otherwise $x_1>x_0$.  If $n_1>m_1$, then
$$
x_1\leq2^{-n_1} + 2^{-n_1-1} + \cdots + 2^{-n_1-a+1}
= 2^{1-n_1} \big( 1-2^{-a} \big)<2^{-m_1}(1-2^{-a})\leq x_0- 2^{-m_k-a},
$$
contrary to hypothesis.  Thus $n_1=m_1$.  If $k=1$, then $l=1$ as $x_1\leq x_0$. Otherwise, we must have $n_2\geq m_2$. Again if $n_2>m_2$, then
\[ x_1  \leq  2^{-n_1} + 2^{-n_2} + 2^{-n_2-1} + \cdots + 2^{-n_2-a+2}
\leq 2^{-m_1} + 2^{-m_2}\big(1-2^{1-a}\big)\leq x_0-2^{-m_k-a}, \]
gives a contradiction. Therefore $n_2=m_2$.  Proceeding in this way shows that $l=k$ and $n_1=m_1,\cdots,n_k=m_k$, establishing the claim.

A similar argument shows that any $x_1\in X_{a,\alpha}$ with $x_0\leq x_1<x_0+2^{-m_k-a}$ can be written in the form
$$
x_1=\sum_{i=1}^k2^{-m_i}+\sum_{i=k+1}^l2^{-m_i}
$$
for some $k\leq l\leq a$ and $m_k<m_{k+1}<\cdots<m_l$.  If in addition $x_1<x_0+2^{-\beta}$, then $m_{k+1}\geq \beta$ so $x_1\in Y$. Thus $(x_0-2^{-m_k-a},x_0+\min(2^{-m_k-a},2^{-\beta}))\cap X_{a,\alpha}=Y$ and so $Y$ is open in $X_{a,\alpha}$.
\end{proof}

Let $S$ be a subsemigroup of a commutative group $(G,+)$, and let $S_1=\mathbb Z_+\times S$.  Let $(w_n)_{n=1}^\infty$ be a sequence in $S$ with the property that for each $s,t\in S$, $s-t\in G$ can be written in at most one way as
\begin{equation}\label{gc:one}
s-t=\sum_{i=1}^k w_{m_i} - \sum_{j=1}^l w_{n_j} \tag{*}\end{equation}
where $k,l\in\mathbb Z_+$ and $\{ m_1,\ldots,m_k, n_1,\ldots,n_l \}$ is a collection of
distinct natural numbers. As usual the empty sum takes the value $0\in G$ so the condition (*) shows that if $k,l\geq 1$ and 
$$
\sum_{i=1}^kw_{m_i}=\sum_{j=1}^lw_{n_i}
$$
for some sets (of distinct natural numbers) $\{m_1,\dots,m_k\}$ and $\{n_1,\dots,n_l\}$, then $k=l$ and $\{m_1,\dots,m_k\}=\{n_1,\dots,n_l\}$. For $(a,s)\in S_1$, and $\alpha\in\mathbb N$, define
\[ U_{a,s,\alpha}=
    \left\{\Big(a-k, s+\sum_{i=1}^k w_{m_i}\Big):\  0\le k\le a,\
      \alpha\le m_1<\cdots<m_k\right\}. \]

\begin{lemma}\label{gen_con_one}
With the notation and conditions as above, the collection $\{U_{a,s,\alpha}\}$ forms a base for a topology $\sigma$ on $S_1$ making $(S_1,\sigma)$ a locally compact semitopological semigroup.
\end{lemma}
\begin{proof}
We first show that $\{U_{a,s,\alpha}:\ (a,s)\in S_1,\ \alpha\in\mathbb N\}$ is a base for a topology $\sigma$ on $S_1$. To this end, note that if $(b,t)\in U_{a,s,\alpha}\setminus\{(a,s)\}$, say $(b,t)=(a-k, s+\sum_{i=1}^k w_{m_i})$,
then $U_{b,t,\beta}\subset U_{a,s, \alpha}$ whenever $\beta\ge m_k+1$. Now for $(a,s), (b,t)\in S_1$ and 
$\alpha,\beta\in\mathbb N$, take $(c,u) \in  U_{a,s,\alpha} \cap U_{b,t,\beta}$.  Then
there exists $\gamma_1, \gamma_2$ such that $U_{c,u,\gamma_1} \subseteq U_{a,s,\alpha}$ and
$U_{c,u,\gamma_2} \subseteq U_{b,t,\beta}$.  Then
\[ U_{c,u,\max(\gamma_1,\gamma_2)} \subseteq U_{c,u,\gamma_1} \cap U_{c,u,\gamma_2}
\subseteq U_{a,s,\alpha} \cap U_{b,t,\beta}, \]
so the $U_{a,s,\alpha}$ do form a base for a topology.

Next we show that $\sigma$ is a Hausdorff topology.  Take $(a,s)\neq (b,t)\in S_1$ and suppose first that $s=t$ so that $a\not=b$.  If there exists $\alpha$ and $\beta$ with $U_{a,s,\alpha}\cap U_{b,t,\beta}\not=\emptyset$, then we can find
$k,l$ with $a-k=b-l$, and sequences $(m_i)$ and $(n_j)$ with
\[ s + \sum_{i=1}^k \omega_{m_i} = s + \sum_{j=1}^l \omega_{n_j}. \]
Condition (\ref{gc:one}) enables us to conclude that $k=l$, so that $a=b$, a contradiction. In particular $U_{a,s,1}$ and $U_{b,s,1}$ are disjoint neighbourhoods of $(a,s)$ and $(b,s)$ respectively. Now suppose that $t-s\neq 0$.  If $t-s$ can be written as
\begin{align}\label{standardexpression}
  t-s=\sum_{i=1}^{\tilde{k}} w_{\tilde{m}_i}- \sum_{j=1}^{\tilde{l}} w_{\tilde{n}_j},
\end{align}
where $\tilde{k}+\tilde{l}\ge 1$ and where $\tilde{m}_i$ and $\tilde{n}_j$ are all distinct
natural numbers, then set
\[ \alpha=\beta=\max\{\tilde{m}_1,\ldots,\tilde{m}_{\tilde{k}},
\tilde{n}_1,\ldots,\tilde{n}_{\tilde{l}}\}+1. \]
Otherwise, set $\alpha=\beta=1$.  Assume toward a contradiction that
$U_{a,s,\alpha}\cap U_{b,t,\beta}\neq\emptyset$.  So we can find $k,l$ with $a-k=b-l$,
and $(m_i)$ and $(n_j)$ with
\begin{align*}
    s+\sum_{i=1}^k w_{m_i}=t+\sum_{j=1}^l
    w_{n_j}\quad \text{that is,} \quad
    t-s=\sum_{i=1}^k w_{m_i}- \sum_{j=1}^l w_{n_j}.
\end{align*}
The uniqueness of the expression (\ref{standardexpression}) implies that
\[ \{\tilde{m}_1,\ldots, \tilde{m}_{\tilde{k}},
\tilde{n}_1,\ldots,\tilde{n}_{\tilde{l}}\}\subset
\{m_1,\ldots, m_k, n_1,\ldots, n_l\}. \]
This contradicts either $\alpha\le m_1< \ldots < m_k$ or $\beta\le n_1<\ldots < n_l$.
In conclusion $\sigma$ is a Hausdorff topology.

Our next objective is to show that $\sigma$ is locally compact.
For $a\in\mathbb Z_+$ and $\alpha\in\mathbb N$ define
\[ X_{a,\alpha} = \{0\}\ \cup\ \Big\{\sum_{i=1}^k 2^{-m_i}:\ 1\le
k\le a,\ \alpha\le m_1<\ldots < m_k \Big\} \subseteq [0,1], \]
so that $X_{a,\alpha}$ is compact.  Fix $(a,s,\alpha)$ and define a map
\[ \psi : U_{a,s,\alpha} \rightarrow X_{a,\alpha}, \quad
(a,s) \mapsto 0, \
\Big(a-k,s+\sum_{i=1}^k w_{m_i}\Big) \mapsto
\sum_{i=1}^k 2^{-m_i}. \]
Condition (\ref{gc:one}) implies that $\psi$ is well-defined; it is then obvious that
$\psi$ is a bijection.  We claim that $\psi$ is actually a homeomorphism.  As
$U_{a,s,\alpha}$ is Hausdorff and $X_{a,\alpha}$ is compact, it is enough to show
that $\psi^{-1}$ is continuous, or equivalently, that $\psi$ is open.  

To show this,
as sets of the form $U_{b,t,\beta}$ form a basis for $\sigma$,
it is enough to show that for each $(b,t)\in U_{a,s,\alpha}$ and $\beta$ with
$U_{b,t,\beta} \subseteq U_{a,s,\alpha}$, we have that $\psi(U_{b,t,\beta})$ is open
in $X_{a,\alpha}$. For $(b,t)\in U_{a,s,\alpha}\setminus\{(a,s)\}$, say $(b,t)=(a-k, s+\sum_{i=1}^k w_{m_i})$ we have $U_{b,t,\beta} \subseteq U_{a,s,\alpha}$ if and only if $\beta>m_k$.  Then
\[ \psi(U_{b,t,\beta}) = \Big\{\sum_{i=1}^{k+l} 2^{-m_i}:\ 0\le l\le a-k,\
\beta\le m_{k+1}<\ldots<m_{k+l}\Big\}, \]
which is open in $X_{a,\alpha}$ by Lemma~\ref{technical}.  Similarly $\psi(U_{a,s,\beta})$ is open in $X_{a,\alpha}$ for all $\beta\geq\alpha$, and so $\sigma$ is locally compact.

Finally we show that $\sigma$ makes the semigroup operation separately continuous. Let $(b,t)\in S$. We claim that the map
\[ M_{b,t}: (x,u)\mapsto (x+b,u+t),\quad S_1\to S_1, \]
is $\sigma$-continuous. Indeed, let $(a,s)\in S_1$ and let $\alpha\in \mathbb N$.
Then $(x,u)\in M_{b,t}^{-1}(U_{a+b,s+t,\alpha})$ if and only if
\[ x+b=a+b-k\,,\quad\textrm{and}\quad u+t=s+t+\sum_{i=1}^k w_{m_i}; \]
where $0\le k\le a+b$,  and $\alpha \le m_1<\ldots< m_k$. This
happens if and only if
\[ x=a-k\,,\quad\textrm{and}\quad u=s+\sum_{i=1}^k w_{m_i}; \]
where $0\le k\le a$,  and $\alpha \le m_1<\ldots< m_k$. So
$M_{b,t}^{-1}(U_{a+b,s+t,\alpha})=U_{a,s,\alpha}$. Since
$(a,s)\in S_1$ and $\alpha\in \mathbb N$ are arbitrary, we deduce that $M_{b,t}$
is continuous.
\end{proof}

We continue with the notation introduced prior to Lemma~\ref{gen_con_one}.  We say (**) holds if for each $t\in S$, there exists an $\alpha_t\in\mathbb N$
such that whenever $s\in S$, $n\ge\alpha_t$, and $s-t+w_n\in S$, then $s-t\in S$.

\begin{lemma}\label{gen_con_two}
With the notation above, if both (\ref{gc:one}) and (**) hold, then the topology $\sigma$ satisfies the property that a sequence $(x_n,u_n)$ in
$S_1$ converges whenever $(x_n+a,u_n+s)$ converges, for some $(a,s)\in S_1$.
\end{lemma}
\begin{proof}
Let $(a,s)\in S_1$ and let $(x_n,u_n)\subseteq S_1$ be such that $(x_n+a,u_n+s)$
converges to some $(b,t)\in S_1$. For each $\alpha\ge\alpha_s\in\mathbb N$, there
exist $n_\alpha$ such that whenever $n\ge n_\alpha$ we have that $(x_n+a,u_n+s)
\in U_{b,t,\alpha}$.  That is, there exists $0\leq l\leq b$ and $\alpha \leq n_1<\cdots
< n_l$ with
\[ x_n+a=b-l,\quad u_n+s=t+\sum_{j=1}^l w_{n_j}. \]

So $b-l\geq a$.  Also, as $n_l \geq \alpha \geq \alpha_s$, and $t-s+\sum_{i=1}^{l-1}
w_{n_i} + w_{n_l} = u_n \in S$, (**) shows that $t-s+\sum_{i=1}^{l-1}
w_{n_i} \in S$.  By induction, $t-s\in S$.  Thus
\[ x_n=b-a-l,\quad u_n=t-s+\sum_{j=1}^l w_{n_j}, \]
where $0\le l\le b-a$ and $\alpha \le n_1<\ldots< n_l$; that is,
$(x_n,u_n)\in U_{b-a,t-s,\alpha}$.  We conclude that $(x_n,u_n)$
converges to $(b-a,t-s)$, as required.
\end{proof}

\begin{remark}\label{rem4}
When $(w_n)$ satisfies the conditions (\ref{gc:one}) and (**), so too does any subsequence of $(w_n)$.  For a subsequence $(w_{m_i})$, denote by
$\sigma_{(m_i)}$ the topology on $S_1$ constructed from $(w_{m_i})$.  We see that
$(0,w_{m_i})\to (1,0)$ with respect to $\sigma_{(m_i)}$.  If $(w_{n_j})$ is another
subsequence such that $\{m_i\}\triangle\{n_j\}$ is infinite, then it is easy to
see that $(0,w_{n_j}) \not\rightarrow (1,0)$ with respect to $\sigma_{(m_i)}$.  Thus the
topologies $\sigma_{(n_j)}$ and $\sigma_{(m_i)}$ differ.  We can then deduce that there
exists a continuum of different topologies $\sigma$ on $S_1$ satisfying the conclusions
of the previous two lemmas.
\end{remark}

Corollary~\ref{cor:one} above shows that $\ell^1(\mathbb Z_+^2)$ has a unique Hopf
algebra predual.  We shall now show that $\mathbb Z_+^2$ admits a continuum of distinct
non-discrete locally compact topologies: of course, none can satisfy the 2nd condition of
Corollary~\ref{good_corr} or Lemma~\ref{best_semigroup}.  We need the following easy fact.
\begin{proposition}
Let $(m_i: 1\le i\le k)$ and $(n_j: 1\le j\le l)$ be sequences in $\mathbb N$ such that
each $k\in\mathbb N$ occurs at most twice in each sequence.  Suppose that
\[ \sum_{i=1}^k 2^{2^{m_i}}=\sum_{j=1}^l 2^{2^{n_j}}. \]
Then $k=l$, and $(m_i)$ and $(n_j)$ are rearrangements of each other.
\end{proposition}

\begin{theorem}\label{infinitegenerators}
There exist a continuum of locally compact topologies on $(\mathbb Z_+^2,+)$
making it a semitopological semigroup.
\end{theorem}
\begin{proof}
Set $G=\mathbb Z$ and $S=\mathbb Z_+$, and let $w_n = 2^{2^n}$ for $n\in\mathbb N$.
The previous proposition shows that this sequence satisfies (*) so
Lemma~\ref{gen_con_one} and Remark~\ref{rem4} apply.
\end{proof}

\begin{theorem}\label{3generators}
There exist a continuum of Hopf algebra preduals for
$\ell^1(\mathbb Z_+\times\mathbb Z)$.
\end{theorem}
\begin{proof}
Now we set $G=S=\mathbb Z$, and have $(w_n)$ as above.  As $S=G$, the condition (**) is obviously satisfied.
The result then follows from Lemma~\ref{best_semigroup}.
\end{proof}

\begin{theorem}\label{Ncdot}
There exist a continuum of Hopf algebra preduals for
$\ell^1(\mathbb N,\cdot)$.
\end{theorem}
\begin{proof}
Let $G$ be the multiplicative group of positive rational numbers, and let $S$ be the
sub-semigroup consisting of odd natural numbers.  We see that $\mathbb Z_+ \times S \cong
(\mathbb N,\cdot)$ by the isomorphism $(k,n) \mapsto 2^k n$.
Finally, let $(w_n)$ be any increasing sequence of odd prime numbers. We see that
$(w_n)$ satisfies the conditions (*) and (**). The result again follows from Lemma~\ref{best_semigroup}.
\end{proof}

\section{Semigroup algebras with no Hopf algebra preduals}

Recall that a semigroup $S$ is said to be an \emph{inverse semigroup} if,
for each $s\in S$, there exists a unique $s^{-1}\in S$ such that $ss^{-1}s=s$
and $s^{-1}ss^{-1}=s^{-1}$.  In this section, we shall exhibit an
inverse semigroup $S$ which admits no semitopological structure giving rise to a
Hopf algebra predual for $\ell^1(S)$.

Let $S$ be the collection of maps $f:\mathbb N_0\rightarrow\mathbb N_0$ such
that there exists a finite (and possibly empty) set $F\subset\mathbb N$ such that $f$ maps $F$ injectively into $\mathbb N$ and $f(n)=0$ for $n\not\in F$. We call $F$ the \emph{injective domain} of $f$. Under composition of functions $S$ is a countable inverse semigroup.

\begin{lemma}
Let $\sigma$ be any locally compact Hausdorff topology on $S$ making
$S$ into a semitopological semigroup.  Then there exists $f\in S$, with injective
domain $F\subset\mathbb N$ such that, for any finite set $F'\subset\mathbb N$
disjoint from $F$, the set
\[ \mc O(f,F') = \{ h\in S : h(n)=f(n) \ (n\in F), \ h(n)=0 \ (n\in F') \}. \]
is open.
\end{lemma}
\begin{proof}
As $S$ is countable, and $(S,\sigma)$ is a Baire space, there exists $f\in S$
with $\{f\}$ open.  As $(S,\sigma)$ is
semitopological, for each $g\in S$, the set $\{h\in S : hg=f \}$ is open.

Fix $F'\subset\mathbb N$ finite and disjoint from $F$.  Define $g\in S$ with injective
domain $G$ so that $g(n)=n$ for $n\in F$ and $g(\mathbb N\cap(G\setminus F))=F'$.  Then
$h\in S$ has $hg=f$ if and only if $h(n)=f(n)$ for $n\in F$ and $h(n)=0$ for $n\in F'$,
that is, if and only if $h\in\mathcal O(f,F')$.
\end{proof}

\begin{theorem}\label{no_hopf_predual}
There is no locally compact Hausdorff topology $\sigma$ on $S$ such that
$C_0(S,\sigma)$ is a predual for $\ell^1(S)$.
\end{theorem}
\begin{proof}
Suppose towards a contradiction that $C_0(S,\sigma)$ is a predual for
$\ell^1(S)$.  By the previous lemma, there exists $f\in S$ with injective domain $F$ such that $\mc O(f,F')$ is open for any finite $F'\subset\mathbb N$ which is disjoint from $F$. Fix such $F'$ and fix $n_0\in F'$.
For each $n$, define $f_n:\mathbb N_0\rightarrow\mathbb N_0$ by
\[ f_n(k) = \begin{cases} f(k) &: k\in F, \\ n &: k=n_0, \\ 0 &:\text{otherwise}.
\end{cases} \]
Then $f_n\in S$ for sufficiently large $n$.  Let $p\in S$ be the
function $p(k)=k$ for $k\in F$, and $p(k)=0$ otherwise.  

Suppose, towards a contradiction, that $(f_n)$ does not converge to $f$.
Then there exists an open set $U\in\sigma$ with $f\in U$, and a sequence
$n_1<n_2<\cdots$ with $f_{n_i} \not\in U$ for all $i$.  Now, $f_np = f$ for
all $n$, so $f_{n_i}p \rightarrow f$.  By Corollary~\ref{good_corr}, there
exists a subsequence of $(f_{n_i})$ which converges to $g$, say.  As $\sigma$
is Hausdorff, $g\neq f$.  However, as $f_{n_i}p=f$ for all $i$, it follows that $gp=f$. That is $g(k)=f(k)$ for all $k\in F$.  As $g\neq f$ there exists $k_0\not\in F$ with $g(k_0)\neq 0$.

Let $h_1\in S$ be the function $h_1(k_0)=k_0$ and $h_1(k)=0$ otherwise, and let $h_2\in S$
be the function $h_2(g(k_0))=g(k_0)$ and $h_2(k)=0$ otherwise.
By construction $h_2 g h_1 \not = 0$, yet $h_2 f_{n_i} h_1 = 0$ unless
$k_0=n_0$ and $n_i = g(n_0)$. In particular $h_2f_{n_i}h_1=0$ for sufficiently large $i$.  Thus
\[ 0 \not= h_2 g h_1 = \lim_i h_2 f_{n_i} h_1 = 0, \]
a contradiction, as required.

Therefore $f_n \rightarrow f$.  However, $f\in \mc O(f,F')$, yet
$f_n \not\in \mc O(f,F')$ for all $n$, giving the required contradiction
to finish the proof.
\end{proof}

\section{A unique algebraic predual}

We end with an example of an infinite semigroup $S$ for which $c_0(S)$ is the unique Banach algebra predual on $\ell^1(S)$.  Other examples of Banach algebras for which the predual is uniquely determined include von Neumann algebras (Sakai's classical result shows that von Neumann algebras have a unique \emph{isometric} predual.  The extension to the non-isometric case can be found in \cite{our}) and $\mathcal B(E)$ for a reflexive Banach space $E$ with the approximation property, \cite[Theorem 4.4]{Daws}.

We first need a little machinery.  Let $\mc A$ be a
Banach algebra.  We turn $\mc A^*$ into an $\mc A$-bimodule in the usual way,
\[ \ip{a\cdot\mu}{b} = \ip{\mu}{ba}, \quad
\ip{\mu\cdot a}{b} = \ip{\mu}{ab} \qquad (a,b\in\mc A, \mu\in\mc A^*). \]
Define bilinear maps from $\mc A^{**}\times\mc A^*$ and $\mc A^*\times
\mc A^{**}$ to $\mc A^*$ by
\[ \ip{\Phi\cdot\mu}{a} = \ip{\Phi}{\mu\cdot a}, \quad
\ip{\mu\cdot\Phi}{a} = \ip{\Phi}{a\cdot\mu}
\qquad (a\in\mc A, \mu\in\mc A^*, \Phi\in\mc A^{**}). \]
Finally, define bilinear maps $\aone,\atwo:\mc A^{**}\times\mc A^{**}
\rightarrow\mc A^{**}$ by
\[ \ip{\Phi\aone\Psi}{\mu} = \ip{\Phi}{\Psi\cdot\mu}, \quad
\ip{\Phi\atwo\Psi}{\mu} = \ip{\Psi}{\mu\cdot\Phi}
\qquad (\Phi,\Psi\in\mc A^{**}, \mu\in\mc A^*). \]
These are the \emph{Arens products}; they are contractive Banach algebra
products on $\mc A^{**}$.  For further details, see \cite[Section~2.6]{dales} or
\cite[Section~1.4]{palmer}.
When $\mc A$ is commutative, $\Phi\aone\Psi =
\Psi\atwo\Phi$ for $\Phi,\Psi\in\mc A^{**}$.

Define $\wap(\mc A^*) \subseteq \mc A^*$
to be those functionals $\mu\in\mc A^*$ such that
\[ \ip{\Phi\aone\Psi}{\mu} = \ip{\Phi\atwo\Psi}{\mu}
\qquad (\Phi,\Psi\in\mc A^{**}). \]
Then $\wap(\mc A^*)$ is an $\mc A$-submodule of $\mc A^*$.
It is a simple calculation (see \cite[Section~2]{Daws} or \cite[Section~4]{Runde2})
that if $\mc A$ is a dual Banach algebra with predual $\mc A_*\subseteq\mc A^*$,
then $\mc A_* \subseteq \wap(\mc A^*)$.

Now let $S = (\mathbb N,\max)$.
Consider the character space of $\ell^\infty(S)$, which as $S$ is discrete,
is equal to the space of ultrafilters on $S$, written $\beta S$.  For $\omega\in
\beta S$, let $\delta_\omega\in\ell^\infty(S)^*$ be the character induced
by $\omega$, so that
\[ \ip{\delta_\omega}{f} = \lim_{s\rightarrow\omega} f(s)
\qquad (f\in\ell^\infty(S)). \]
So for $S=(\mathbb N,\max)$, let $\omega,\upsilon\in\beta S$ be non-principal,
so that for $f\in\ell^\infty(S)$, we have that
\[ \ip{\delta_\omega\aone\delta_\upsilon}{f}
= \lim_{s\rightarrow\omega} \lim_{t\rightarrow\upsilon} f(\max(s,t))
= \lim_{t\rightarrow\upsilon} f(t) = \ip{\delta_\upsilon}{f}. \]
Hence, if $f\in\wap(\ell^\infty(S))\subseteq\ell^\infty(S)$, then
\[ \ip{\delta_\upsilon}{f} = \ip{\delta_\omega\aone\delta_\upsilon}{f}
= \ip{\delta_\upsilon\aone\delta_\omega}{f} = \ip{\delta_\omega}{f}. \]
It follows easily that
\[ \wap(\ell^\infty(\mathbb N,\max)) = c_0(\mathbb N) \oplus \mathbb C1. \]

\begin{lemma}
Let $\mc A_*$ be a Banach space, let $\mc A = \mc A_*^*$, and let
$F\subseteq \mc A^*$ be a closed subspace such that $\mc A_*\subseteq F$ and
$F/\mc A_*$ is one-dimensional.  Let $E\subseteq F$ be a
spacial predual for $\mc A$.  Then
\[ E^\perp := \big\{ M\in F^* : \ip{M}{\mu}=0 \ (\mu\in E) \big\} \]
is also one-dimensional.
\end{lemma}
\begin{proof}
We can find $\mu_0\in\mc A^*\setminus\mc A_*$ with $F$ being the span of $\mc A_*$
and $\mu_0$.  Pick $M_0\in\mc A^{**}$ with $\ip{M_0}{\mu_0}=1$ and $\ip{M_0}{\mu}=0$
for all $\mu\in\mc A_*$.  By restriction, we shall regard $M_0$ as a member of $F^*$.
If $E^\perp=\{0\}$ then $E=F$, which is a contradiction,
as $F$ strictly contains $\mc A_*$ and so cannot be a predual for $\mc A$.
So, towards a contradiction, suppose that we can find linearly independent vectors
$M_1,M_2\in E^\perp$.

For $i=1,2$, if we restrict $M_i$ to $\mc A_* \subseteq F$, then we induce a member
of $\mc A_*^* = \mc A$, say $a_i\in\mc A$, which satisfies $\ip{M_i}{\mu} = \ip{\mu}{a_i}$ for
$\mu\in\mc A_*$.  Then $M_i - a_i$ annihilates $\mc A_*$, so as $F$ is the linear span
of $\mc A_*$ and $\mu_0$, we can find $\alpha_i\in\mathbb C$ with $M_i - a_i = \alpha_i M_0$.
We hence have that
\[ 0 = \ip{M_i}{\mu} = \ip{\mu}{a_i} + \alpha_i \ip{M_0}{\mu}
\qquad ( \mu\in E, i=1,2 ). \]
If $\alpha_1=0$, then for each $\mu\in E$, we have that $\ip{\mu}{a_1} = 0$.
As $E$ is a predual for $\mc A$, this means that $a_1=0$, so that $M_1=0$, a contradiction.
Similarly, $\alpha_2\not=0$.

We hence see that
\[ \ip{\mu}{\alpha_1^{-1} a_1} = - \ip{M_0}{\mu} = \ip{\mu}{\alpha_2^{-1} a_2}
\qquad (\mu\in E). \]
As $E$ is a predual, this shows that $\alpha_1^{-1}a_1 = \alpha_1^{-2}a_2$.  Thus
\begin{align*} M_1 &= a_1 + \alpha_1 M_0 = \alpha_1\big( \alpha_1^{-1}a_1 + M_0 \big)
= \alpha_1\big( \alpha_2^{-1} a_2 + M_0 \big) \\
&= \alpha_1\alpha_2^{-1} \big( a_2 + \alpha_2 M_0 \big)
= \alpha_1\alpha_2^{-1} M_2, \end{align*}
a contradiction, as required.
\end{proof}

\begin{theorem}\label{Unique.l1s}
Let $S=(\mathbb N,\max)$, and let $E\subseteq\ell^\infty(S)$ be a
predual for $\ell^1(S)$.  Then $E = c_0(S)$.
\end{theorem}
\begin{proof}
We have that
\[ E \subseteq \wap(\ell^\infty(S)) = c_0(\mathbb N) \oplus \mathbb C1. \]
We identify the dual of $c_0(\mathbb N) \oplus \mathbb C 1$ with $\ell^1(\mathbb N)
\oplus \mathbb C1$.  By the previous lemma, $E^\perp$ is one dimensional, so there
exists $a\in\ell^1(\mathbb N)$ and $\alpha\in\mathbb C$, not both zero, such that
\[ E^\perp = \{ \Phi\in \wap(\ell^\infty(S))^* : \ip{\Phi}{\mu}=0 \ (\mu\in E) \}
= \mathbb C( a+\alpha 1). \]
It hence follows that
\[ E = \big\{ (x,\beta)\in c_0 \oplus \mathbb C1 : \ip{a}{x} = -\alpha\beta \big\}. \]
Then $E=c_0$ if and only if $a=0$.

So, towards a contradiction, suppose that $a\not=0$.  Pick $(x,\beta)\in E$.
As $1\in\ell^\infty(S)$ is clearly invariant for the $\ell^1(S)$ module action, we see
that $(\delta_s\cdot x,\beta)\in E$ for $s\in S$, and so
\[ \ip{a}{x} = -\alpha\beta = \ip{a}{\delta_s\cdot x} = \ip{\delta_s\cdot a}{x}
\qquad (s\in S). \]
Let $a = \sum_{n\in\mathbb N} a_n \delta_n$, so for $s\in\mathbb N$,
\[ \delta_s \cdot a = \sum_n a_n \delta_{\max(s,n)}
= \Big(\sum_{n=1}^s a_n\Big) \delta_s + \sum_{n=s+1}^\infty a_n \delta_n. \]
Let $x = \sum_n x_n \delta_n$, so we see that
\[ \sum_{n=1}^\infty a_n x_n = \ip{a}{x} = \ip{\delta_s\cdot a}{x}
= \sum_{n=1}^s a_n x_s + \sum_{n=s+1}^\infty a_n x_n
\qquad (s\in\mathbb N). \]
That is,
\[ \sum_{n=1}^s a_n x_n = x_s \sum_{n=1}^s a_n \qquad (s\in\mathbb N). \]
Letting $s\rightarrow\infty$, we conclude that $\ip{a}{x} = \sum_n a_n x_n = 0$.
As $(x,\beta)\in E$, we see that $\ip{a}{x} = -\alpha\beta$, so either $\alpha=0$,
or $\beta=0$.  If $(x,\beta)\in E$ implies that $\beta=0$, then $E\subseteq c_0$,
which as $E$ is a predual means that $E=c_0$ as required.

Otherwise, we have that $\alpha=0$, so that $(x,\beta)\in E$ if and only if
$\ip{a}{x}=0$.  If $\ip{1}{a} = 0$, then for $(x,\beta)\in E$, $\ip{x+\beta 1}{a} = 0$,
so $a$ annihilates $E$.  As $E$ is a predual, $a=0$, contradiction.  So $\sum_n a_n\not=0$.
As $E$ is an $\ell^1(S)$-module, we have that $\ip{a}{x}=0$ implies
that $\ip{\delta_s\cdot a}{x}=0$.  If $a = a_{s_0} \delta_{s_0}$ for some $s_0\in S$, then
$\ip{x}{a}=0$ if and only if $\ip{x}{\delta_{s_0}}=0$, which clearly does not imply
that $\ip{x}{\delta_s\cdot\delta_{s_0}}$ is zero for all $s$.
Otherwise, choose $s_0<s_1$ minimal with $a_{s_0} \not=0$ and $a_{s_1}\not=0$.  Let $s$
be greater than $s_0$ and $s_1$ chosen such that $\sum_{n=1}^s a_n \not=0$, which is
possible, as $\sum_n a_n \not=0$.  Let
$x = a_{s_0} \delta_{s_1} - a_{s_1}\delta_{s_0} + \delta_s$, so that $\ip{x}{a}=0$, but
\[ \ip{x}{\delta_s\cdot a} = \sum_{n=1}^s a_n \not = 0. \]
This final contradiction completes the proof.
\end{proof}

The underlying fact which allows this proof to work is that for $S=(\mathbb N,\max)$,
we have that $\wap(\ell^\infty(S))$ is very small.  In \cite{Chou}, Chou shows that
when $G$ is an infinite discrete group, then $\wap(\ell^\infty(G)) / c_0(G)$ contains
an isometric copy of $\ell^\infty$.  So there is no hope of a generalisation of
the above proof to group algebras.

\noindent Matthew Daws  \\
School of Mathematics  \\
University of Leeds  \\
Leeds LS2 9JT  \\
United Kingdom  \\
Email: \texttt{matt.daws@cantab.net}

\medskip

\noindent Hung Le Pham  \\
Department of Mathematical and Statistical Sciences  \\
University of Alberta  \\
Edmonton T6G 2E1  \\
Canada  \\
Email: \texttt{hlpham@math.ualberta.ca}

\medskip

\noindent Stuart White  \\
Department of Mathematics  \\
University of Glasgow  \\
Glasgow G12 8QW  \\
United Kingdom  \\
Email: \texttt{s.white@maths.gla.ac.uk}

\end{document}